\documentclass[reqno,10pt]{amsart}

\usepackage{amssymb, amsmath, amsthm}
\usepackage{color}
 \usepackage[bookmarks=false]{hyperref}
\hypersetup{
    colorlinks=true, 
    linktoc=all,     
    linkcolor=blue,  
    citecolor=blue
}

\newcommand*{\fplus}{\genfrac{}{}{0pt}{}{}{+}}
\newcommand*{\fdots}{\genfrac{}{}{0pt}{}{}{\cdots}}
\newcommand*{\fminus}{\genfrac{}{}{0pt}{}{}{-}}

\newcommand{\qrfac}[2]{{\left({#1}; q\right)_{#2}}} 
\newcommand{\pqrfac}[3]{{\left({#1};#3\right)_{#2}}}

\newcommand{\qbin}{\genfrac{[}{]}{0pt}{}}

\newmuskip\pFqskip
\mathchardef\pFcomma=\mathcode`, 

\newmuskip\pFqskip
\mathchardef\pFcomma=\mathcode`, 

\newtheorem*{rems}{Remarks} 
\newenvironment{Remarks}{\begin{rems}\normalfont}{\end{rems}}

\newtheorem*{rem}{Remark} 
\newenvironment{Remark}{\begin{rem}\normalfont}{\end{rem}}

\numberwithin{equation}{section}

\allowdisplaybreaks

\dedicatory{Dedicated to the memory of Dick Askey}

\begin{document} 
\title[Ramanujan's $q$-continued fractions]{Ramanujan's $q$-continued fractions}

\author[G.~Bhatnagar]{Gaurav Bhatnagar
}
\address{Ashoka University, Sonipat, Haryana 131029, India}
\email{bhatnagarg@gmail.com}
\urladdr{https://www.gbhatnagar.com}
%


\date{\today}

\begin{abstract}
Ramanujan's $q$-continued fractions are a central part of Ramanujan's development of basic hypergeometric series.
They appear in Chapter 16 of Part III and Chapter 32 of Part V of {\em Ramanujan's Notebooks} edited by Berndt, and in Volume I of Andrews and Berndt's {\em Ramanujan's Lost Notebook}. In these references the continued fractions as presented in the order in which they appear in Ramanujan's original notebooks. 
We  summarize the work of several authors on this topic and re-organize Ramanujan's $q$-continued fractions. 
\end{abstract}


\maketitle

\section{Introduction}

Ramanujan considered the $q$-analogue of the well-known Fibonacci continued fraction
\begin{equation*}\label{golden-mean-cfrac}
\frac{1}{1}\fplus \frac{1}{1}\fplus\frac{1}{1}\fplus\fdots
\end{equation*}
given by
\begin{subequations}
\begin{equation}\label{RRcfrac1}
\frac{1}{1}\fplus\frac{q}{1}\fplus\frac{q^2}{1}\fplus\frac{q^3}{1}\fplus\fdots,
\end{equation}
and showed that for $|q|<1$, this continued fraction can be written as a ratio of two very similar sums:
\begin{equation}
 \frac
{\displaystyle\sum_{k=0}^{\infty} \frac{q^{k^2+k}}{{(1-q)(1-q^2)\cdots (1-q^k)}}}
{\displaystyle\sum_{k=0}^{\infty} \frac{q^{k^2}}{(1-q)(1-q^2)\cdots (1-q^k)}}.
 \label{rr-cfrac1}
\end{equation}
The continued fraction \eqref{RRcfrac1} is called the Rogers--Ramanujan continued fraction. In view of the Rogers--Ramanujan identities, this ratio can be expressed as a ratio of infinite products
\begin{equation}
\frac
{(1-q)(1-q^6)(1-q^{11}) \cdots  }
{(1-q^2)(1-q^7)(1-q^{12}) \cdots  }
\times
\frac{(1-q^4)(1-q^9) (1-q^{14}) \cdots}
{(1-q^3) (1-q^8)(1-q^{13}) \cdots}
.\label{entry16.38.iii}
\end{equation}
\end{subequations}

It is clear that Ramanujan considered his work on this continued fraction as one of the highlights of his work. 
Consider the last two entries of Chapter 16 of Ramanujan's second notebook \cite{RamanujanNB}. Entry 38  concerns the continued fraction above and the Rogers--Ramanujan identities. Entry 39 computes special cases (see Berndt \cite[p.~77--86]{Berndt1991}). As a corollary, Ramanujan obtains
\begin{subequations}
\begin{align}
 \frac{1}{1}\fplus\frac{e^{-2\pi}}{1}\fplus\frac{e^{-4\pi}}{1}
\fplus\frac{e^{-6\pi}}{1}\fplus\fdots
   &=\Bigg( \sqrt{\frac{5+\sqrt{5}}{2} } - \frac{\sqrt{5}+1}{2} \Bigg) e^{2\pi/5},
\label{II.16.39ii-cor}
\\ 
\intertext{and}
 \frac{1}{1}\fminus\frac{e^{-\pi}}{1}\fplus\frac{e^{-2\pi}}{1}
\fminus\frac{e^{-3\pi}}{1}\fplus\fdots
 &= \Bigg(\sqrt{\frac{5- \sqrt{5}}{2} } - \frac{\sqrt{5}-1}{2}\Bigg) e^{\pi/5}  .
\label{II.16.39i-cor}
\end{align}
\end{subequations}

Askey \cite{Askey1989} suggested that Ramanujan's $q$-extension 
\eqref{RRcfrac1} was the primary motivation for Ramanujan to discover the Rogers--Ramanujan identities, and indeed for Ramanujan's development of basic hypergeometric series. In support of this hypothesis, we note that
Ramanujan ended his development of basic hypergeometric series with 
\eqref{II.16.39ii-cor} and \eqref{II.16.39i-cor}. Thus $q$-continued fractions may be quite central to Ramanujan's thought process.  

Ramanujan himself included these two continued fractions   in his first letter \cite[p.~29]{BR1995} to Hardy in 1913; this is additional evidence that these were results Ramanujan liked and thought likely to impress Hardy. 

And impress him, they did. About these, and another formula related to this continued fraction, Hardy \cite[p.~9]{Hardy1959} wrote
\begin{quote}
(they) defeated me completely; I had never seen anything in the least like them before. A single look at them is enough to show that they could only be written by a mathematician of the highest class.
\end{quote}

The objective of this article is to summarize the work of many mathematicians, whose work has made it easier to understand and organize Ramanujan's continued fractions. %
 The three formulas \eqref{RRcfrac1}--\eqref{entry16.38.iii} serve as a model for this organization.
 

As it turns out, there are only five general results. In all of these a ratio of two series is expressed as a continued fraction. The remaining continued fractions are special cases of these. These general results are given in \S\S\ref{g-sec}, \ref{G-sec} and \ref{sec:entry11-12}. In some cases, the same ratio of two series is
expressed as two or three continued fractions; thus, they are transformation formulas too. We note some special cases in \S\ref{sec:cf-trans}. There are some special cases where the ratio reduces to a single series, and many where they 
reduce to a ratio of infinite products (see \S\ref{sec:cf-prods}). When this happens, such as in the case of the Rogers--Ramanujan 
continued fraction, Ramanujan computes special values.  
In this article, we list all of Ramanujan's $q$-continued fractions. 
( Formulas for the computation of special values have not been considered here.) 

 
Apart from this organization, we are interested in the ideas involved in the proofs. A variety of clever ideas are required to prove  Ramanujan's $q$-continued fractions; we present a sample, 
giving details when possible, and references otherwise.
We conclude in \S\ref{sec:proofs} with some remarks arising from our study, and pointers to recent work in this area.

\subsection*{Notation}
The {\em $q$-rising factorial}  is defined as $\qrfac{a}{0} :=1$, 
$$\qrfac{a}{n} := 
(1-a)(1-aq)\cdots (1-aq^{n-1})  \text{ for } n=1, 2, \dots;
$$ 
and, for $|q|<1$,
$$\qrfac{a}{\infty} := \prod_{k=0}^\infty (1-aq^{k}).$$
In addition, we use the short-hand notation
\begin{align*}
\qrfac{a_1, a_2,\dots, a_r}{k} &:= \qrfac{a_1}{k} \qrfac{a_2}{k}\cdots
\qrfac{a_r}{k}.
\end{align*}
We will require the $q$-binomial coefficient, defined as
$$\qbin{n}{k}_q := \frac{\qrfac{q}{n}}{\qrfac{q}{k}\qrfac{q}{n-k}},$$
where $n\geq k$ are non-negative integers. When $n<k$, we take $\qbin{n}{k}_q=0.$

\subsection*{Ramanujan's entries}
We record the notation used to refer to Ramanujan's entries from \cite{RamanujanNB, Ramanujan-LN}. An entry takes its number from either one of the five volumes of Ramanujan's Notebooks edited by Berndt, or the Lost Notebook edited by Andrews and Berndt. We have prefixed the volume number and/or the chapter to the entry. 
For example, 
Entry III.16.39 refers to Entry 39, of Chapter 16 in Part III of Berndt's five volumes on {\em Ramanujan's notebooks}; that is, it refers to \cite[Ch.~16, Entry 39]{Berndt1991}. 
We place an additional L when referring to the Lost Notebook. 
Thus, L.I.6.2.1 refers to Entry 6.2.1 in Part I of the {\em Ramanujan's lost notebook} edited by Andrews and Berndt, that is,  \cite[Entry 6.2.1]{AB2005}. 

\section{A more general continued fraction}\label{g-sec}

In his second letter to Hardy, Ramanujan~\cite[p.~57]{BR1995} mentions that Entry~III.16.29~(i)
is ``a particular case of a theorem on
the continued fraction 
$$\frac{1}{1}\fplus \frac{ax}{1}\fplus\frac{ax^2}{1}\fplus\frac{ax^3}{1}\fplus\frac{ax^4}{1}\fplus
\frac{ax^5}{1}\fplus
\fdots ,$$
which is a particular case of the continued fraction
$$
\frac{1}{1}\fplus\frac{a x}{1+bx}\fplus\frac{a x^2}{1+bx^2}\fplus\frac{a x^3}{1+bx^3}\fplus
 \fdots, 
$$
which is a particular case of a general theorem on continued fractions.'' 
A theorem regarding this continued fraction involves the sum
\begin{equation}
g(b,\lambda):=\sum_{k=0}^{\infty}\frac{\lambda^k q^{k^2}}{\qrfac{q,-bq}{k}}.\label{g}
\end{equation}
For the expressions corresponding to \eqref{g} in the letter, Ramanujan had $x$ in place of $q$, and $\lambda$
in place of $a$.  

Ramanujan records the following continued fractions in the Lost Notebook \cite[p.~40]{Ramanujan-LN}
\subsection*{L.I.6.3.1 (i), (ii), (iii)} Let $|q|<1$. Then
\begin{subequations}
\begin{align}
\frac{g(b,\lambda q)}{g(b,\lambda)}
&
=\frac{1}{1}\fplus\frac{\lambda q}{1}\fplus\frac{bq+\lambda q^2 }{1}\fplus\frac{\lambda q^3}{1}\fplus
\frac{bq^2+\lambda q^4}{1}\fplus\fdots  \label{g-cfrac1}
\\
&=\frac{1}{1}\fplus\frac{\lambda q}{1+bq}\fplus\frac{\lambda q^2}{1+bq^2}\fplus\frac{\lambda q^3}{1+bq^3}\fplus
\frac{\lambda q^4}{1+bq^4}\fplus \fdots \label{g-cfrac2}.\\
\intertext{Let $|q|<1$ and $|b|< 1$. Then}
\frac{g(b,\lambda q)}{g(b,\lambda)}
&=
\frac{1}{1-b}\fplus\frac{b+\lambda q}{1-b}\fplus\frac{b+\lambda q^2}{1-b}\fplus\frac{b+\lambda q^3}{1-b}\fplus
\frac{b+\lambda q^4}{1-b}\fplus \fdots.
\label{g-cfrac3} 
\end{align}
\end{subequations}
\begin{Remarks}
A continued fraction of Touchard~\cite{Touchard1952} can be considered to be a special case of \eqref{g-cfrac3}, see Prodinger~\cite{Prodinger2012}.
The first two continued fractions can be obtained as special cases of the general continued fractions in \S\ref{G-sec}. We prove the third. 
\end{Remarks}
\begin{proof}[Proof of \eqref{g-cfrac3}]  
We present J.~Cigler's proof (sent by email, Sept.~2012),
which illustrates an idea of Askey~\cite{Askey1989} 
regarding a discovery approach to the Rogers--Ramanujan identities beginning with \eqref{RRcfrac1} (see also \cite{GB2015}). 

Let $F(b,\lambda)$ denote the continued fraction \eqref{g-cfrac3}. Then 
$$F(b,\lambda)=\frac{1}{1-b+(b+\lambda q)F(b,\lambda q)}.$$
To change this non-linear equation into a linear one, we substitute
$$F(b,\lambda)=\frac{g(b,\lambda q)}{g(b,\lambda )},$$
to obtain Ramanujan's recurrence relation given in 
Entry L.I.6.3.1(iv):
\begin{equation*}
g(b,\lambda) = (1-b)g(b,\lambda q)  +(b+\lambda q) g(b,\lambda q^2).
\end{equation*}
Now assume that 
$$g(b,\lambda)=\sum_{n=0}^\infty a_n\lambda^n,$$
substitute in Ramanujan's recurrence relation, and compare coefficients of $\lambda^n$, to obtain
$$a_n = (1-b)q^n a_n + bq^{2n} a_n +q^{2n-1}a_{n-1},$$
or
$$
a_n = \frac{q^{2n-1}}{(1-q^n)(1+bq^n)} a_{n-1}.
$$
Iterating this we obtain
$$ a_n = \frac{q^{n^2}}{\pqrfac{q, -bq}{n}{q}} a_0.$$
Now taking $a_0=1$, we obtain the expression \eqref{g} for $g(b,\lambda)$ and a (formal) proof of
\eqref{g-cfrac3}.
\end{proof}

The continued fraction  \eqref{g-cfrac2} appeared earlier as Entry III.16.15.
When $b=0$  
this reduces to a relation for an extension of the Rogers--Ramanujan continued fraction.
\subsection*{Cor. to III.16.15} If $|q|<1$, then
\begin{equation*}\label{cor-entry15}
\frac{\displaystyle\sum_{k=0}^{\infty} \frac{\lambda^k q^{k^2+k}}{\qrfac{q}{k}}}
{\displaystyle\sum_{k=0}^{\infty} \frac{\lambda^k q^{k^2}}{\qrfac{q}{k}}}
=
\frac{1}{1}\fplus \frac{\lambda q}{1}\fplus\frac{\lambda q^2}{1}\fplus\frac{\lambda q^3}{1}
\fplus\frac{\lambda q^4}{1}\fplus
\fdots .
\end{equation*}
When $\lambda=1$, this reduces to the Rogers--Ramanujan continued fraction. In the very next entry Ramanujan gives a formula for the convergents of the (reciprocal of) this continued fraction.
\subsection*{III.16.16}
\begin{equation*}
\frac{\mu_n(0)}{\mu_n(1)} =
1+ \frac{\lambda q}{1}\fplus
\frac{\lambda q^2}{1}\fplus\frac{\lambda q^3}{1}
\fplus\fdots \fplus \frac{\lambda q^n}{1}
 ,\label{entry16}
\end{equation*}
where
$\mu_n(s)$ is defined as
\begin{equation*}\label{mu1-finite}
\mu_n(s):=
\sum_{k=0}^{\infty} q^{k^2+sk}\lambda^k
\qbin{n-k-s+1}{k}_{q}
, 
\end{equation*}
 for $s=0, 1, 2, \dots.$
 \begin{Remarks}\
In our proof we illustrate a method due to Euler~\cite{LE1788-616} (see \cite{GB2014, BI2021}), that works on all the general $q$-continued fractions mentioned in this article. 
Euler's approach is based on using the elementary identity
\begin{equation*}\label{div-1step}
\frac{N}{D}=1+\frac{N-D}{D}
\end{equation*}
to \lq divide' two series, both of whose first terms are $1$.
\end{Remarks}

%
\begin{proof} [Proof of III.16.16]
Note that the sum $\mu_n(s)$ is finite; the index $k$ goes from $0$ to $ \lfloor \frac{n-s+1}{2}\rfloor$.
We will show
\begin{equation}\label{gn-recursion}
\frac{\mu_n(s)}{\mu_n(s+1)}=
1+\frac{\lambda  q^{s+1}}{{\cfrac{\mu_n(s+1)}{\mu_n(s+2)}}},
\end{equation}
for $s=0,1, 2, 3, \dots, n-1$.
Formula \eqref{entry16} follows by taking $s=0, 1, \dots, n-1$ in turn. Observe that the iteration stops when $s=n-1$, because
\begin{equation*}
\mu_n(n)=1=\mu_n(n+1).
\end{equation*} 

To prove \eqref{gn-recursion}, note that
\begin{align*}
\frac{\mu_n(s)}{\mu_n(s+1)} &=1+\frac{\mu_n(s)-\mu_n(s+1)}{\mu_n(s+1)} \\
& = 1+\frac{1}{\mu_n(s+1)} \\
&\hspace{0.35cm}\times
\displaystyle\sum_{k=0}^{\infty}  
\frac{q^{k^2+sk}\lambda^k\qrfac{q}{n-k-s}}{\qrfac{q}{k}\qrfac{q}{n-2k-s+1}}
  \left[(1-q^{n-k-s+1})- (1-q^{n-2k-s+1})q^k\right]
 \\
&= 
1+ \frac{1}{\mu_n(s+1)}
\displaystyle\sum_{k=1}^{\infty}  
\frac{q^{k^2+sk}\lambda^k}{\qrfac{q}{k-1}}
\frac{\qrfac{q}{n-k-s}}{\qrfac{q}{n-2k-s+1}}
\\
&= 
1+\frac{\lambda q^{s+1}}{\mu_n(s+1)}
\displaystyle\sum_{k=0}^{\infty}  
\frac{q^{k^2+(s+2)k)}\lambda^{k}}{\qrfac{q}{k}}
\frac{\qrfac{q}{n-k-s-1}}{\qrfac{q}{n-2k-s-1}},
\end{align*}
on shifting the sum. From this \eqref{gn-recursion} follows and the proof is complete. 
\end{proof}


\subsection*{III.16.13 and Corollary 6.2.4 to L.I.6.2.1(Eisenstein)}If $|q|<1$, then
\begin{equation*}
\sum_{k=0}^{\infty} (-a)^k q^{\frac{k(k+1)}{2}}=
\frac{1}{1}\fplus\frac{aq}{1}\fplus\frac{a(q^2 - q)}{1}\fplus\frac{aq^3}{1}\fplus\frac{a(q^4-q^2)}{1}
\fplus\fdots ,\label{entry13}
\end{equation*}
\begin{Remarks} Eisenstein's continued fraction appears in \cite[pp.~35--39]{Eisenstein1975}.  The case $a=1$ was known to Gauss in 1797, see \cite[p.~152]{AB2005}. 
See Folsom~\cite{Folsom2006} for a history,  and an account of Eisenstein's continued fractions.
\end{Remarks}
\begin{proof}[Proof outline] 
The special  case $\lambda=a$ and $b=-a$ of \eqref{g-cfrac1} can be written as Entry III.16.13, 
where to obtain the sum on the left hand side, we have to apply the transformation formula III.16.9 to the sums. 
\end{proof}
Ramanujan (see \cite{Berndt1991}) noted formulas for the denominators of the convergents of this continued fraction. They are given by the following: 
\begin{align*}
D_{2n} &= \sum_{k=0}^n a^{k} q^{nk} \qbin{n}{k}_q, \\
\intertext{and}
D_{2n+1} &= \sum_{k=0}^n a^{k} q^{(n+1)k} \qbin{n}{k}_q .
\end{align*}

%
%

\section{Still more general}\label{G-sec}
Some candidates for what Ramanujan may have meant by the phrase `general  theorem on continued fractions' in his letter, concern the sum 
\begin{equation}\label{G}
G(a,b,\lambda):=\sum_{k=0}^{\infty}
\frac{\qrfac{-\lambda/a}{k} \, a^k\, q^{(k^2+k)/2}}{\qrfac{q, -bq}{k}}.
\end{equation}
\begin{subequations}
However, only special cases of these theorems appear in Ramanujan's notebooks, which is the only record of 
his discoveries before he wrote the aforementioned letter to Hardy. 
In the lost notebook, Ramanujan gave the following continued fractions and derived numerous special cases.
\subsection*{L.I.6.2.1 and L.I.6.4.1} Let $a$, $b$, $\lambda$, and $q$ be complex numbers, with $|q|<1$. Then
\begin{align}\label{G-cfrac1}
&\hspace{-1cm}\frac{G(aq, b, \lambda q)}{G(a, b, \lambda)}\cr
&=\frac{1}{1}\fplus\frac{aq+\lambda q}{1}\fplus  \frac{bq+\lambda q^2 }{1}\fplus\frac{aq^2+\lambda q^3}{1}\fplus\frac{bq^2+\lambda q^4}{1} \fdots \\
&= \frac{1}{1+aq}\fplus\frac{\lambda q - abq^2}{1+aq^2+bq}\fplus\frac{\lambda q^2 -abq^4 }
{1+aq^3+bq^2}\fplus\frac{\lambda q^3-abq^6}{1+aq^4+bq^3}\fplus \fdots . \label{G-cfrac2}
\end{align}
\begin{Remarks}
\ 
\begin{enumerate}
\item The continued fraction in \eqref{G-cfrac2}, Entry L.I.6.4.1, appears a few pages after 
\eqref{G-cfrac1} in the Lost Notebook \cite{Ramanujan-LN}, as a transformation formula given by the second equality above. 
\item Note that if we take $a=0$ in \eqref{G-cfrac1} we obtain  \eqref{g-cfrac1}; and taking $a=0$ in \eqref{G-cfrac2}
we obtain \eqref{g-cfrac2}. But there is no general continued fraction in Ramanujan's notebooks for 
${G(aq, b, \lambda q)}/{G(a, b, \lambda)}$ which reduces to \eqref{g-cfrac3}. 
A candidate for such a continued fraction is ~\cite[Theorem~6.4.1]{AB2005}
\begin{multline}\label{G-cfrac3} 
\frac{G(aq, b, \lambda q)}{G(a, b, \lambda)}
=\\
\frac{1}{1}\fplus\frac{aq+\lambda q}{1-aq+bq}\fplus\frac{aq+\lambda q^2}
{1-aq+bq^2}\fplus
\frac{aq+ \lambda q^3}{1-aq+bq^3}\fplus \fdots, \\
 \text{ where } |aq|<1  .
\end{multline} 
When $b=0$ and $a\mapsto b/q$, this reduces to \eqref{g-cfrac3}. 
This continued fraction appears in different forms in Hirschhorn \cite{MDH1980} and  Bhargava and Adiga \cite{BA1984}. We need to  appeal to Entry III.16.8 in order to match the two forms. 
\item In addition to \eqref{G-cfrac1}, \eqref{G-cfrac2}, and \eqref{G-cfrac3}, there are two more continued fractions for ${G(aq, b, \lambda q)}/{G(a, b, \lambda)}$ in \cite{GB2014} (see also \cite{LMS2017}).
\item Another candidate for a `general theorem on continued fractions' is 
Andrews \cite[Theorem~6]{Andrews1968}. The most general continued fraction of this type is due to
Andrews and Bowman~\cite{AndBow1995}; see also Gupta and Masson~\cite{GM1999}.
\end{enumerate}
\end{Remarks}
\end{subequations}

\begin{proof}[On proofs of \eqref{G-cfrac1}--\eqref{G-cfrac3}]
Andrews \cite{andrews1979} included a proof of \eqref{G-cfrac1} in his introductory article reporting his discovery of the Lost Notebook. Bhargava and Adiga~\cite{BA1984} proved \eqref{G-cfrac2} and 
\eqref{G-cfrac3}.
\end{proof}
There are several special cases of these continued fractions noted by Ramanujan. In each of the following, we assume the conditions of L.I.6.2.1 (if applicable).

\subsection*{L.I.6.2.3}
\begin{equation*}\label{G-cfrac1-cor}
\frac{G(a, b, 0)}{G(aq, b, 0)}
=1+\frac{aq}{1}\fplus\frac{bq }{1}\fplus\frac{aq^2}{1}\fplus\frac{bq^2}{1} \fdots.
\end{equation*}
\begin{proof}
Set $\lambda =0$ in \eqref{G-cfrac1} and take reciprocals.
\end{proof}
\subsection*{Cor. 6.2.9 to L.I.6.2.1}
\begin{equation*}
\sum_{k=0}^\infty (-1)^k q^{3k^2+2k} \big(1+q^{2k+1}\big)
=
\frac{1}{1}\fplus\frac{q^2-q}{1}\fplus\frac{ q^4 - q^2}{1}\fplus\frac{q^6-q^3}{1}\fplus
\frac{ q^8 - q^4}{1}\fplus\fdots .
\label{LNIcor6.2.9} 
\end{equation*}
\begin{proof}
We take $q\mapsto q^2$ and $a=-1/q$, $b=-1$, and $\lambda =1$ in \eqref{G-cfrac1}.
The sum side requires an iterate of Heine's transformation formula, and L.I.9.5.1,
see \cite[p.~155]{AB2005}. Heine's formula is Entry III.16.6.
\end{proof}

\subsection*{Cor.~6.2.11 to Entry L.I.6.2.1}
\begin{equation*}
1-\sum_{k=1}^\infty q^{k(3k-1)/2} \big(1-q^{k}\big)
=
\frac{2}{2}\fplus\frac{q+q}{1}\fplus\frac{ q^2 + q^3}{1}\fplus\frac{q^3+q^5}{1}\fplus
\frac{ q^4 + q^7}{1}\fplus\fdots .
\label{LNIcor6.2.11} 
\end{equation*}
\begin{proof} Take $q\mapsto q^2$, $a= 1/q$, $b=1$, and $\lambda = 1/q,$ 
in \eqref{G-cfrac1}. The sum side requires elementary series manipulations, a special case of the $q$-binomial theorem (III.16.2), the second iterate of Heine's transformation (\cite[eq.~(III.3)]{GR90}) and Entry L.I.9.4.7. See \cite[p.~156]{AB2005} for the details. 
\end{proof}
There are further special cases in the next two sections.

\section{Transformations of continued fractions} \label{sec:cf-trans}
In \S\S\ref{g-sec} and \ref{G-sec}, we saw that a ratio of two series  can be written in multiple ways as continued fractions. 
In this section we note two  transformation formulas recorded by Ramanujan, obtained by setting the same special case of two of these continued fractions equal to one another. In addition, there is one transformation which follows by taking the odd part of a continued fraction. 

\subsection*{L.I.6.5.1}
Let $k\geq 0$, $\alpha = (1+\sqrt{1+4k} )/2$, and $\beta = (-1+\sqrt{1+4k} )/2$. Then for $|q|<1$,
\begin{subequations}
\begin{align}
\frac{1}{1}\fplus & \frac{k+q}{1}\fplus  \frac{k+q^2}{1}\fplus \frac{k+q^3}{1}\fplus\fdots
\label{LI.6.5.1a} \\
&=
\frac{1}{\alpha}\fplus \frac{q}{\alpha+\beta q}\fplus \frac{q^2}{\alpha+\beta q^2}\fplus \frac{q^3}{\alpha+\beta q^3}
\fplus\fdots .\label{LI.6.5.1b} 
\end{align}
\end{subequations}
\begin{proof}
Observe that $\alpha\beta=k$ and $\alpha-\beta =1$. To show this transformation formula we take a special case $\lambda = 1/\alpha^2$, $b=\beta/\alpha$ in  \eqref{g-cfrac2} and \eqref{g-cfrac3} and equate the two resulting continued fractions.  
After  substituting in \eqref{g-cfrac2}, we multiply and divide each fraction by $\alpha$, and simplify to obtain
$\alpha$ times \eqref{LI.6.5.1b}. Similarly, we obtain $\alpha$ times \eqref{LI.6.5.1a} from \eqref{g-cfrac3}.
\end{proof}

%
%
%
 
\subsection*{L.I.6.5.2} For $|q|<1$,
\begin{align*}
\frac{1}{1}\fplus & \frac{2+q}{1}\fplus  \frac{2+q^2}{1}\fplus \frac{2+q^3}{1}\fplus\fdots \cr
&=
\frac{1}{2}\fplus \frac{q}{2+ q}\fplus \frac{q^2}{2+q^2}\fplus \frac{q^3}{2+ q^3}
\fplus\fdots.
\end{align*}
\begin{proof}
This is a special case of Entry L.I.6.5.1, where $k=2$, so that $\alpha=2$ and $\beta =1$.
\end{proof}

\subsection*{L.I.6.4.2} If $|q|<1$, then
\begin{align*}
\frac{1}{a+c}\fminus & \frac{ab}{a+b+cq}\fminus\fdots\fminus \frac{ab}{a+b+cq^n}\fminus\fdots 
\\
&= \frac{1}{c-b+a}\fplus  \frac{bc}{c-b+a/q}\fplus\fdots\fplus \frac{bc}{c-b+a/q^n}\fplus\fdots 
\end{align*}
\begin{proof} Take $\lambda=0$, $a\mapsto -b/aq$, $b\mapsto c/a$ in \eqref{G-cfrac3} and \eqref{G-cfrac2}, and compare the resulting continued fractions. The transformation formula is obtained after some simplification, see \cite[p.~161]{AB2005} for more details. 
\end{proof}

\subsection*{L.I.6.4.3} For $q\neq 0$, $n$ a positive integer, we have
\begin{align}\label{LI6.4.3}
1+\frac{a}{1}\fplus & \frac{b}{q}\fplus\frac{a}{1}\fplus \frac{b}{q^2}\fdots\fplus \frac{b}{q^n}\fplus\frac{a}{1} 
\cr
&=
1+a - \frac{ab}{a+b+q}\fminus \frac{ab}{a+b+q^2}\fminus\fdots\fminus \frac{ab}{a+b+q^n}.
\end{align}
\begin{proof} The proof is elementary. The continued fraction on the left hand side of  \eqref{LI6.4.3} is of the form
$1+{N_{2n+1}}/{D_{2n+1}}$, where
$$\frac{N_{2n+1}}{D_{2n+1}}
= \frac{a_1}{b_1}\fplus\frac{a_2}{b_2}\fplus\fdots \fplus \frac{a_{2n+1}}{b_{2n+1}}.$$
We consider instead the 
continued fraction with convergents $N_k/D_k$.
The continued fraction on the right hand side is obtained by adding $1$ to the continued fraction with 
convergents $N_1/D_1$, $N_3/D_3$, $\dots$, ${N_{2n+1}}/{D_{2n+1}}$.
This is known as the odd part of the continued fraction.

The details are as follows. An equivalent form of Ramanujan's Entry II.12.1 is that $N_k$ and $D_k$ are both determined by a three-term recurrence relation of the form
\begin{align*}
P_k &= b_kP_{k-1} +a_kP_{k-2}, 
\end{align*}
for $k = 3, 4, \dots, 2n+1$, with initial values $N_1=a_1$, $D_1=b_1$, $N_2=a_1b_1$ and $D_2=a_2+b_1b_2$. 
It is easy to see that $N_{2k+1}$ and $D_{2k+1}$ (for $k=2, 3, \dots, n$) also satisfy the same recurrence relation
with $a_k \mapsto - {a_{2k-1}a_{2k}b_{2k+1}}/{b_{2k-1}}$, and 
$b_k \mapsto a_{2k+1}+b_{2k}b_{2k+1}+{a_{2k}b_{2k+1}}/{b_{2k-1}}.$
Now we substitute the values of $a_k$ and $b_k$ from the left hand side of \eqref{LI6.4.3}, and find the corresponding recurrence for the odd part of the continued fraction. On adding $1$ to the corresponding continued fraction, we obtain the 
right hand side of \eqref{LI6.4.3}.
\end{proof}


\section{Continued fractions written as products}\label{sec:cf-prods}
So far we have seen continued fractions written as a ratio of two series. Many of Ramanujan's $q$-continued fractions are ratios of infinite products. It turns out that all of these are obtained when the series in question can be written as infinite products. 
An example is the Rogers--Ramanujan continued fraction highlighted in the introduction, which
has the three expressions given by formulas \eqref{RRcfrac1},
\eqref{rr-cfrac1} and \eqref{entry16.38.iii}.
In this section, we list such special cases of the general continued fractions given so far. In all the formulas in this section, we assume $|q|<1$.

\subsection*{Corollary~6.2.2 to L.I.6.2.1}
For any complex number $a\neq -q^{2n+1}$, $n\ge 1$
\begin{equation}
\frac{\pqrfac{-aq^2}{\infty}{q^2}}{(-aq; q^2)^{ 2}_{\infty}}
=
\frac{1}{1}\fplus\frac{aq}{1}\fplus\frac{ q + aq^2}{1}\fplus\frac{aq^3}{1}\fplus
\frac{ q^2 + aq^4}{1}\fplus\fdots .
\label{LNIcor6.2.2} 
\end{equation}
\begin{proof} Take $a=0$, $b=1$ and $\lambda\mapsto a$ in \eqref{G-cfrac1}. 
This requires the $q$-binomial theorem (Entry III.16.2) to sum both the sums. 
\end{proof}
The next two continued fractions are further special cases of this continued fraction.
\subsection*{V.32.21; Corollary~6.2.1 to L.I.6.2.1}
\begin{equation*}
\frac{\pqrfac{-q^2}{\infty}{q^2}}{(-q; q^2)_{\infty}} =
\frac{\pqrfac{q}{\infty}{q^2}}{(q^2; q^4)^{ 2}_{\infty}}
=
\frac{1}{1}\fplus\frac{q}{1}\fplus\frac{ q + q^2}{1}\fplus\frac{q^3}{1}\fplus
\frac{ q^2 + q^4}{1}\fplus\fdots .
\label{LNIcor6.2.1} 
\end{equation*}
\begin{proof}
Take $a = 1$ in \eqref{LNIcor6.2.2}.
\end{proof}

\subsection*{Corollary~6.2.10 to L.I.6.2.1}
\begin{equation}
\frac{\pqrfac{-q^3}{\infty}{q^4}}{(-q; q^4)_{\infty}}
=
\frac{1}{1}\fplus\frac{q}{1}\fplus\frac{ q^2 + q^3}{1}\fplus\frac{q^5}{1}\fplus
\frac{ q^4 + q^7}{1}\fplus\fdots .
\label{LNIcor6.2.10} 
\end{equation}
\begin{proof}
Take $q\mapsto q^2$ and $a= 1/q$ in \eqref{LNIcor6.2.2}.
\end{proof}

Another continued fraction for the LHS of \eqref{LNIcor6.2.10} is as follows.
\subsection*{V.32.20}
\begin{equation*}
\frac{\pqrfac{q^3}{\infty}{q^4}}{(q; q^4)_{\infty}}
=
\frac{1}{1}\fminus\frac{q}{1+q^2}\fminus\frac{ q^3}{1+q^4}\fminus\frac{q^5}{1+q^6}\fminus
\frac{q^7}{1+q^8}\fminus
\fdots .
\end{equation*}
\begin{proof}
This is obtained from \eqref{g-cfrac2} by taking $q\mapsto q^2$, $\lambda = -1/q$, $b=1$. The sums are special case of the $q$-binomial theorem (the case $a\mapsto 0$ of Entry III.16.2, with $q\mapsto q^4$).
\end{proof}

Next, we have two more special cases of \eqref{G-cfrac1}. 
\subsection*{Cor.~6.2.7 to L.I.6.2.1; V.32.18}
\begin{equation*}
\frac{\pqrfac{q, q^5}{\infty}{q^6}}{(q^3; q^6)^2_{\infty}}
= \frac{\pqrfac{q}{\infty}{q^2}}{(q^3; q^6)^3_{\infty}}
=
\frac{1}{1}\fplus\frac{q+q^2}{1}\fplus\frac{ q^2 + q^4}{1}\fplus\frac{q^3+q^6}{1}\fplus
\frac{ q^4 + q^8}{1}\fplus\fdots .
\label{LNIcor6.2.7} 
\end{equation*}
\begin{Remark} This continued fraction is known as Ramanujan's cubic continued fraction. See 
Chan \cite{Chan1995} and \cite[Chapter 3]{AB2005}. 
\end{Remark}
\begin{proof}
Take $q\mapsto q^2$, and take $a=1/q$, $b=1$, and $\lambda =1$ in \eqref{G-cfrac1}.
To obtain the product side, we require some identities of Slater \cite{Slater1951, Slater1952}. For the details,
see \cite[p.~154]{AB2005}. 
\end{proof}

\subsection*{V.32.22; V.32.23; Cor.~6.2.8 to L.I.6.2.1}
\begin{equation*}
\frac{\pqrfac{q, q^7}{\infty}{q^8}}{(q^3, q^5; q^8)_{\infty}}
=
\frac{1}{1}\fplus\frac{q+q^2}{1}\fplus\frac{ q^4}{1}\fplus\frac{q^3+q^6}{1}\fplus
\frac{ q^8}{1}\fplus\fdots .
\label{LNIcor6.2.8} 
\end{equation*}
\begin{proof}
Take $q\mapsto q^2$, and take $a=1/q$, $b=0$, and $\lambda =1$ in \eqref{G-cfrac1}. To obtain the product side, we require some identities of Slater \cite{Slater1952}. 
See \cite[p.~154]{AB2005} for more details. This continued fraction is referred to as the Ramanujan--G{\" o}llnitz--Gordon continued fraction. 
\end{proof}

\subsection*{V.32.19}
\begin{equation}
\frac{\pqrfac{q^2}{\infty}{q^3}}{\pqrfac{q}{\infty}{q^3}}
=
\frac{1}{1}\fminus\frac{q}{1+q}\fminus\frac{ q^3}{1+q^2}\fminus\frac{q^5}{1+q^3}\fminus
\frac{ q^7}{1+q^4}\fminus\fdots .
\label{V32.19} 
\end{equation}
\begin{Remark} The product side of \eqref{V32.19} and a related continued fraction appeared in an intriguing claim of Ramanujan concerning a continued fraction that converges to three different limits. This claim and related ideas have been explored by Andrews, Berndt, Sohn, Yee and Zaharescu~\cite{ABSYZ2003, ABSYZ2005} and Ismail and Stanton \cite{IS2006}.
\end{Remark}
\begin{proof}
The continued fraction in this entry can be found from \eqref{G-cfrac2}. 
Let $\omega = e^{2\pi i /3}$. Then take 
$\lambda =0$, $a=-\omega/q$, $b=-\omega^2$, to obtain the continued fraction
$$\frac{1}{1-\omega}\fminus\frac{q}{1+q} \fminus\frac{ q^3}{1+q^2}\fminus\fdots,$$
from which \eqref{V32.19} can be obtained by taking reciprocals and adding $\omega$.
A proof along these lines can be found in Berndt~\cite[p.~46]{Berndt1998}. See also \cite{ABSYZ2003}. See the recurrence relation for the function 
$G_2(s)$  in \cite[p.~64]{GB2014} for the connection with \eqref{G-cfrac2}. 
\end{proof}

Note that the proofs of the last three entries require results that are not explicitly stated in Ramanujan's 
notebooks. In fact, these are the only three entries in this article with this property.

\section{Entries III.16.11 and III.16.12}\label{sec:entry11-12}
The first two $q$-continued fractions in Chapter 16 of Ramanujan's second notebook are Entries III.16.11 and III.16.12. Entry III.16.11 is motivated by a continued fraction recorded in Chapter 12 of Ramanujan's second notebook (see also \cite[Chapter~14, Vol.~1]{RamanujanNB}).  In Entry II.12.18, Ramanujan finds a continued fraction for 
$$\frac{(x+1)^n - (x-1)^n}{(x+1)^n +(x-1)^n}.$$
Ramanujan notes four corollaries: continued fractions for $\tan^{-1} x$, $\log\frac{1+x}{1-x}$, $\tan x$, $ \frac{e^x-1}{e^x+1}$, so it would be natural to attempt an extension of II.14.18.

We can expand $(x+1)^n$ and $(x-1)^n$ using the binomial theorem. If $n$ is a non-negative integer, then the numerator becomes a polynomial with terms of the form $a_kx^k$, where $k=n-1, n-3, \dots$;
and the denominator has terms with $k=n, n-2, n-4, \dots$. 

Now the $q$-analog of the binomial theorem, given in Entry III.16.2, is
\begin{equation*}\label{III.16.2}
F(a,b):= \sum_{k=0}^\infty \frac{\qrfac{b/a}{k}}{\qrfac{q}{k}}a^k =
\frac{\qrfac{b}{\infty}}{\qrfac{a}{\infty}}.
\end{equation*}
If we consider 
$$\frac{F(a,b)-F(-a,-b)}{F(a,b)+F(-a,-b)}$$ we obtain
\begin{align*}
 \frac{\displaystyle
\sum_{k=0}^{\infty}\frac{\qrfac{b/a}{2k+1}}
{\qrfac{q}{2k+1}}a^{2k+1}}
{\displaystyle
\sum_{k=0}^{\infty}\frac{\qrfac{b/a}{2k}}
{\qrfac{q}{2k}}a^{2k}}
&=
\frac{{\qrfac{-a, b}{\infty}}-{\qrfac{a, -b}{\infty}}}
{{\qrfac{-a,b}{\infty}}+\qrfac{a, -b}{\infty}}.
\end{align*}
The numerator of the left hand side contains  the odd terms of the $q$-binomial sum, and the denominator contains the even terms. The right hand appears in the following entry.
\subsection*{III.16.11} Let $|q|<1$ and $|a|<1$. Then we have
\begin{multline*}
\frac{{\qrfac{-a, b}{\infty}}-{\qrfac{a, -b}{\infty}}}
{{\qrfac{-a, b}{\infty}}+\qrfac{a, -b}{\infty}}
=
\frac{a-b}{1-q}\fplus\frac{(a-bq)(aq-b)}{1-q^3}\fplus \cr
\frac{q(a-bq^2)(aq^2-b)}{1-q^5}\fplus 
\frac{q(a-bq^3)(aq^3-b)}{1-q^7}\fplus
\fdots .
\end{multline*}
\begin{proof}
For a proof, see Berndt~\cite{Berndt1991}, which relies on ideas of Jacobsen  \cite{LJ1989} to clarify an earlier proof by Adiga, Berndt, Bhargava and Watson~\cite{ABBW1985}. 
\end{proof}
Entry III.16.12 has a very similar continued fraction.  
\subsection*{III.16.12}
Let $|q|<1,$ and $|ab|<1$. Then we have
\begin{multline*}
\frac{\pqrfac{a^2q^3, b^2q^3}{\infty}{q^4}}
{\pqrfac{a^2q, b^2q}{\infty}{q^4}}
=
\frac{1}{1-ab}\fplus\frac{(a-bq)(b-aq)}{(1-ab)(1+q^2)}\fplus \cr
\frac{(a-bq^3)(b-aq^3)}{(1-ab)(1+q^4)}\fplus
\frac{(a-bq^5)(b-aq^5)}{(1-ab)(1+q^6)}\fplus
\fdots.
\end{multline*}
\begin{Remarks}\ 
\begin{enumerate}
\item Entry III.16.12 is a $q$-analogue of Entry II.12.25.
To see this, one has to use III.16.1(ii), which contains the definition of the $q$-gamma function. 
\item The continued fraction in Entry III.16.12 converges for $|ab|>1$. 
It also converges for $|ab|<1$ and $|q|>1$.  The expression for both of these can be found using III.16.12, see
\cite{Berndt1991, BI2021}.
\end{enumerate}
\end{Remarks}
\begin{proof}
For two inter-related proofs, see the work of Ismail and the author  \cite{BI2021}. One of these proofs uses Euler's approach; and the other, a standard method from the theory of orthogonal polynomials. 
\end{proof}
%

This continued fraction also arises from a ratio of two sums.
From the $q$-binomial theorem, it follows that for
 $|q|<1$ and $|a|<1$:
\begin{align*}
\frac{\pqrfac{a^2q^3, b^2q^3}{\infty}{q^4}}
{\pqrfac{a^2q, b^2q}{\infty}{q^4}}
&= 
\frac{\displaystyle
\sum_{k=0}^\infty \frac{\pqrfac{(bq/a)^2}{k}{q^4} }{\pqrfac{q^4}{k}{q^4}} (a^2q)^k
}
{\displaystyle
\sum_{k=0}^\infty \frac{\pqrfac{(b/aq)^2}{k}{q^4} }{\pqrfac{q^4}{k}{q^4}} (a^2q^3)^k
}.
\end{align*}

\section{The proofs of Ramanujan's $q$-continued fractions}\label{sec:proofs}
The first proofs of many of Ramanujan's general $q$-continued fractions were given by  Andrews \cite{andrews1979, An1981} and Adiga, Berndt, Bhargava and Watson \cite{ABBW1985}; proofs are collated in  \cite{AB2005} and \cite{Berndt1991}. Even earlier, Selberg~\cite{Selberg1936} obtained many Rogers--Ramanujan type continued fractions, where the continued fractions are expressed as infinite products. Reuter~\cite{Reuter2014} has surveyed proofs of Ramanujan's hypergeometric continued fractions. 

As we have seen, an eclectic collection of techniques is used to prove Ramanujan's $q$-continued fractions. 
As Berndt~\cite{Berndt2010} wrote in a similar context:
\begin{quote}
Methods for proving these continued fraction formulas are varied and at times ad hoc. Ramanujan evidently had a systematic procedure for proving these continued fraction formulas, but we don't know what it is.
\end{quote}
Nevertheless, when it comes to $q$-continued fractions, there are several interesting systematic approaches. As we have seen, all of Ramanujan's general continued fractions (in \S\S\ref{g-sec}, \ref{G-sec}, \ref{sec:entry11-12}), can be obtained as a ratio of two similar series; the remaining results are obtained as special cases.  In fact, Ramanathan~\cite{KGR1987a} derived all of these general continued fractions formally by finding three-term recurrence relations from contiguous relations for $_2\phi_1$ series.
 
There is, however, a gap in this formal approach to continued fractions. One does not know that the continued fraction converges; and, even if it does, whether it converges to the ratio of series in question. Jacobsen~\cite{LJ1989} has clarified such issues (see \cite{LW1992} for a comprehensive exposition). 

Another way convergence has been addressed is by associating the continued fraction with orthogonal 
polynomials, see Ismail and Stanton~\cite{IS2006}.  However, we require the parameters to be real (rather 
than complex) numbers in this approach.  Continued fractions are closely related to the study of orthogonal 
polynomials (see, for example, Askey and Ismail~\cite{AI1984}). Al--Salam and Ismail \cite{Al-I1983} have 
studied orthogonal polynomials related to Ramanujan's continued fractions, see also \cite{BI2022, IS1997}.  
A third approach to show convergence is to use a theorem of Pincherle
 \cite[p.~164]{JT1980}, see, for examples, Masson~\cite{Masson1989} 
 and Gupta and Masson~\cite{GM1999}.   With regard to questions of convergence of Ramanujan's continued fractions, see also \cite{ABSYZ2003, ABSYZ2005, GB2014, BM2004, BM2007, IS2006}. 

In addition to Ramanathan~\cite{KGR1987a}, Lorentzen (n{\' e}e Jacobsen)~\cite{LL2008}  suggested 
another systematic approach to Ramanujan's  (hypergeometric) continued fractions.  Yet another
approach, following Euler's elementary idea works on all of Ramanujan's general $q$-continued fractions, see~\cite{GB2014, BI2021}. Sokal~\cite{Sokal2022} has developed Euler's idea further and given several examples. 
Hirsch\-horn~\cite{MDH1972, MDH1974, MDH1980} proves many of Ramanujan's continued fractions by finding analogous formulas for the convergents of continued fractions (see also
Menon~\cite{Menon1965}, and \cite{BI2022, BMW2006, Prodinger2018}.)

Regarding transformation formulas for continued fractions, Andrews and Berndt \cite{AB2005} used the 
Bauer--Muir transformation~\cite[p.~76--80]{LW1992} to prove Entry L.I.6.5.1. Lee, Mc Laughlin and Sohn~\cite{LMS2017} have applied this idea to Ramanujan's transformation formulas. By contrast, 
\cite{GB2014} explains the transformation of continued fractions \eqref{G-cfrac2} by transforming the series appearing on the left-hand side of
 \eqref{G-cfrac1}, and then obtaining the corresponding continued fractions by Euler's method. 

As we have seen, taking special cases of Ramanujan's general continued fractions involve the application of 
other transformation and summation formulas. In particular, in \S\ref{sec:cf-prods} we saw that the
Rogers--Ramanujan type identities are very useful in expressing continued fractions as infinite products. 
The first paper with  continued fractions of this kind was written by Selberg~\cite{Selberg1936}. Further examples appear in Gu and Prodinger~\cite{GP2011}. 

This prolific set of ideas surrounding Ramanujan's continued fractions is only to be expected. After all, as Hardy~\cite[p.~XXX]{Ramanujan-CW} famously remarked
\begin{quote}
(Ramanujan's) mastery of continued fractions was, on the formal side at any rate, beyond that of any mathematician in the world...
\end{quote}

\subsection*{Acknowledgements} A shorter version of this article appears in the Encyclopedia of Srinivasa Ramanujan and his Mathematics. We thank Krishnaswami Alladi, George Andrews, Bruce Berndt, Peter Paule, Ole Warnaar, and Ae Ja Yee for their helpful comments.

%



\end{document}